\documentclass{article}
\usepackage{a4wide}
\usepackage{amsfonts}
\usepackage{amsmath}
\usepackage{amsthm}
\usepackage[english]{babel}
\usepackage{bbm}
\usepackage{graphicx}
\usepackage{xcolor}
\usepackage{caption}
\usepackage{subcaption}

\usepackage{algorithm}
\usepackage{algpseudocode}

\newcommand{\E}{\mathbb E}
\newcommand{\N}{\mathbb N}
\renewcommand{\P}{\mathbb P}
\newcommand{\Q}{\mathbb Q}
\newcommand{\R}{\mathbb R}

\newcommand{\half}{\mbox{$\frac 1 2$}}

\newcommand{\1}{\mathbbm 1}

\newcommand{\erf}{\operatorname{erf}}
\newcommand{\erfc}{\operatorname{erfc}}

\newtheorem{theorem}{Theorem}[section]
\newtheorem{lemma}[theorem]{Lemma}
\newtheorem{proposition}[theorem]{Proposition}

\newtheorem{corollary}[theorem]{Corollary}

\theoremstyle{definition}
\newtheorem{definition}[theorem]{Definition}

\newtheorem{problem}[theorem]{Problem}

\theoremstyle{remark}
\newtheorem{remark}[theorem]{Remark}
\newtheorem{example}[theorem]{Example}

\begin{document}
 
\title{Explicit solution of relative entropy weighted control}
\author{Joris Bierkens\thanks{Donders Institute for Brain, Cognition and Behaviour, Radboud University Nijmegen, The Netherlands, j.bierkens@science.ru.nl}, Hilbert J. Kappen\thanks{Donders Institute for Brain, Cognition and Behaviour, Radboud University Nijmegen, The Netherlands, b.kappen@science.ru.nl}}

\maketitle

\begin{abstract}
We consider the minimization over probability measures of the expected value of a random variable, regularized by relative entropy with respect to a given probability distribution. In the general setting we provide a complete characterization of the situations in which a finite optimal value exists and the situations in which a minimizing probability distribution exists. Specializing to the case where the underlying probability distribution is Wiener measure, we characterize finite relative entropy changes of measure in terms of square integrability of the corresponding change of drift. For the optimal change of measure for the relative entropy weighted optimization, an expression involving the Malliavin derivative of the cost random variable is derived. The theory is illustrated by its application to several examples, including the case where the cost variable is the maximum of a standard Brownian motion over a finite time horizon. For this example we obtain an exact optimal drift, as well as an approximation of the optimal drift through a Monte-Carlo algorithm.
\end{abstract}


\thanks{MSC2010 subject classifications: primary 93E20; secondary 60H07, 94A17.} 

\thanks{The research leading to these results has received funding from the European Community's Seventh Framework Programme (FP7/2007-2013) under grant agreement no. 270327}

\section{Introduction}

In certain situations in stochastic optimal control theory, the dynamic programming or Hamilton-Jacobi-Bellman equations may be transformed, through the Hopf-transform, into linear equations \cite{Fleming1982}, \cite[Chapter VI]{FlemingSoner2009}. In the past years, within the applied control and machine learning community, there has been a significant amount of interest in this class of `path integral control problems' (see e.g. \cite{Kappen2005, Todorov2006, Banisch2013}). This class of problems also occurs in risk sensitive control theory (see \cite{FlemingSoner2009}) and the theory of large deviations (see \cite{BoueDupuis1998}), and it occurs in modified form (constrained to equivalent martingale measures) in mathematical finance, in particular as the dual problem for a portfolio optimization problem \cite{Monoyios2013}.  It is the goal of this paper to review and extend the mathematical underpinning of this optimization problem, as well as showcase some new results within this context.

The problem we consider is a minimization problem over probability measures that are absolutely continuous with respect to a given probability measure (referred to as the `uncontrolled measure'). The functional we wish to minimize is the sum of (i) the expectation of a given random variable with respect to any probability measure, and (ii) the relative entropy of that probability measure with respect to the uncontrolled measure. The density of the optimal probability measure with respect to the uncontrolled distribution is readily available through an explicit expression in terms of the cost random variable. The challenge is then to understand this probability measure within the context of the underlying problem.

In particular, in the special case in which we are interested in this paper, the uncontrolled distribution is Wiener measure on the space of continuous sample paths. Absolutely continuous change of distribution then corresponds, by the Girsanov theorem, to a change of drift, which we will interpret as control process. The regularizing relative entropy corresponds to squared control cost. The questions we wish to answer in this paper are: (i) under what conditions does there exist an optimal change of drift corresponding to a given cost functional and probability measure, and (ii) how can it be computed?

There is a close relation to existing theory within the field of large deviations theory and stochastic optimal control \cite{BoueDupuis1998, DupuisEllis1997}. We should also mention the work of F\"ollmer \cite{Follmer1985, Follmer1986}. For an excellent self-contained review of these results, see \cite{Lehec2011}. 
The main aim of this paper is to review the application of the mentioned results in a mathematical control context.  We also illustrate the use of Monte-Carlo methods; for recent work on the use of Monte-Carlo methods in relative entropy weighted control, see \cite{ThijssenKappen2014}.

Furthermore, for the reader who is familiar with the available literature, theoretical contributions of our paper include:
\begin{itemize}
 \item[(i)] Relaxing the usual boundedness assumptions to a condition that guarantees finite relative entropy of the optimal change of measure (condition~\eqref{hyp:finite_entropy} of Section~\ref{sec:RE_optimization});
\item[(ii)] Solution of the problem where the cost random variable is the maximum of a standard Brownian motion with controlled drift over a finite time horizon.
\end{itemize}

\subsection{Outline}
In Section~\ref{sec:RE_optimization}, we consider the general relative entropy weighted optimization problem, and completely characterize the different situations that may arise. The situation with finite relative entropy is most useful, 
and problems for which the optimal change of measure has finite relative entropy are easily characterized in terms of conditions on the cost functional and the probability measure.
Then in Section~\ref{sec:dynamic_RE_optimization}, the finite relative entropy case is further investigated within the context of a Wiener process. It is shown that a change of measure with finite relative entropy corresponds to a square integrable drift, which in particular is the case for the optimal density.
In Section~\ref{sec:control_as_malliavin_derivative} we show how the optimal drift may be computed through the Clark-Ocone formula.
To illustrate the use of this approach, and as an interesting result in its own right, we compute the optimal drift for the case where the cost functional is the maximum of a one dimensional Wiener process with controlled drift on a finite time horizon (Section~\ref{sec:minimize_maximum}). We also provide a Monte-Carlo algorithm for the approximation of such a solution, which is easily extended to other problems. 

\subsection{Notation}
As is common in probability theory, we will allow random variables to assume their values within the extended reals $[-\infty, \infty]$. Resulting formal expressions may be interpreted as follows: $\log 0 = -\infty$, $\log \infty = \infty$, $\exp(-\infty) = 0$, $\exp(\infty) = \infty$ and $\infty \exp(-\infty) = 0$.
For any $a \in \R$, we write $(a)^+ := a \vee 0$ and $(a)^- := -a \vee 0$ for the positive and negative parts of $a$.
The euclidean norm of $x \in \R^d$ is denoted by $|x|$. For an adapted process $\theta$ and a  continuous local martingale $M$, both with values in $\R^d$, we write $\int_0^t \langle \theta_s, d M_s \rangle$ to indicate $\sum_{i=1}^d \int_0^t \theta_s^i \ d M_s^i$.
If $(\Omega,\mathcal F,\P)$ is a probability space we write $\E^{\P}$ for expectation with respect to the probability measure $\P$.
Lebesgue measure will be denoted by $\mathrm{Leb}$.

%

\section{Relative entropy weighted optimization}
\label{sec:RE_optimization}

Let $(\Omega, \mathcal F, \P)$ be a probability space. The probability measure $\P$ will be referred to as the \emph{uncontrolled (probability) measure}. Let $C$ be a random variable assuming values in $[-\infty, \infty]$. The random variable $C$ indicates a cost we wish to minimize, as explained below. 

Let $\mathcal P$ denote the set of probability measures on $(\Omega, \mathcal F)$. We wish to find a probability measure $\Q \in \mathcal P$ that 
\begin{itemize}
\item[(i)] is absolutely continuous with respect to $\P$ (denoted by $\Q \ll \P$), 
\item[(ii)] reduces the expected cost $\E^{\Q} C$, but
\item[(iii)] has small deviation from $\P$. We take the relative entropy
\[ \mathcal H(\Q ; \P) = \int_{\Omega} \log \left( \frac{d \Q}{d \P} \right) \ d \Q = \E^{\Q} \left[ \log \left(\frac{d \Q}{d \P}\right)\right]\] 
as a measure of this deviation (see e.g. \cite[Section 1.4]{DupuisEllis1997} for preliminary results on relative entropy). Recall that $\mathcal H(\Q;\P) \geq 0$ for any $\Q, \P \in \mathcal P$, and $\mathcal H(\Q;\P) = 0$ if and only if $\Q = \P$.
\end{itemize} 
Note that (i) is a constraint and (ii) and (iii) are conflicting optimization targets.

Let $\mathcal P_0 := \left\{ \Q \in \mathcal P : \mbox{$\E^{\Q} [(C)^+] < \infty$ and $\mathcal H(\Q;\P) < \infty$}\right\}$ denote the set of admissable probability measures, and note that $\mathcal P_0$ is convex.
Define the cost functional $J : \mathcal P \rightarrow \R$ by 
\begin{equation} \label{eq:functional} J(\Q) := \left\{ \begin{array}{ll} \E^{\Q} C + \mathcal H(\Q ; \P) = \E^{\Q} \left[ C + \log \left(\frac{d \Q}{d \P}\right)\right] \quad & \mbox{if} \ \Q \in \mathcal P_0, \\
                                                         \infty \quad & \mbox{otherwise}.
                                                        \end{array} \right.
\end{equation}
The definition of $\mathcal P_0$ is as non-restrictive as possible, so that for $\Q \in \mathcal P_0$  the value $J(\Q)$ is well defined within the interval $[-\infty, \infty)$.

We arrive at the following problem:

\begin{problem}[Relative entropy weighted optimization]
\label{prob:REoptimization}
Compute $J^{\star} = \inf_{\Q \in \mathcal P_0} J(\Q)$, and if it exists, a minimizer $\Q^{\star} \in \mathcal P_0$ such that $J(\Q^{\star}) = J^{\star}$.
\end{problem}

The solution of this problem is well known for the case in which $\P(|C| < K) = 1$ for some $K > 0$, 
see e.g. \cite[Proposition 1.4.2]{DupuisEllis1997} or \cite[Proposition 2.5]{BoueDupuis1998}. 
The purpose of this section is to provide a complete characterization of the existence of solutions of Problem~\ref{prob:REoptimization} in terms of conditions on $\P$ and $C$. 
To achieve this goal, we will consider the following further conditions on $C$ and $\P$. 
\begin{equation}
\tag{FE}\label{hyp:finite_entropy} \mbox{\emph{finite (relative) entropy}: } \ \P(C < +\infty) > 0 \quad \mbox{and} \quad \E^{\P}\left[ \exp(-C) |C|  \right] <  \infty.
\end{equation}
\begin{equation}
 \tag{I}\label{hyp:integrability} \mbox{\emph{integrability}: } \ 0 < \E^{\P} \left[ \exp(-C) \right] < \infty.
\end{equation}
Condition~\eqref{hyp:finite_entropy} will earn its name (`finite relative entropy') below: as we will see, it is a necessary and sufficient condition for the `optimal' probability distribution to have finite relative entropy with respect to $\P$.
The implication
\[ \mbox{\eqref{hyp:finite_entropy}} \Longrightarrow \mbox{\eqref{hyp:integrability}}\]
holds as a result of the estimates $\exp(-x) \1_{\{x \leq -1\}} \leq \exp(-x)|x|$ and $\exp(-x) \1_{\{x > -1\}} \leq \frac 1 e$.

\begin{example}
The following examples may serve to illustrate the conditions~\eqref{hyp:finite_entropy} and~\eqref{hyp:integrability}.
\begin{itemize}
 \item[(i)] $\Omega = [0,\infty)$, with $\P$ having density $f(\omega) = \exp(-\omega)$ with respect to Lebesgue measure; $C(\omega) = -\omega$: $\P(C < +\infty) > 0$ but~\eqref{hyp:integrability} does not hold.
\item[(ii)] $\Omega = (-\infty,\infty)$, $\frac{d\P}{d \mathrm{Leb}}(\omega) =\frac{\exp(-|\omega|)}{k(1 + \omega^2)}$, with $k$ a normalization constant, and $C(\omega) = -|\omega|$. Then $\P$ is a probability distribution,
\[ \E^{\P} [\exp(-C)] = \frac 1 k \int_{-\infty}^{\infty} \frac 1 {1 + x^2} \ d x = \frac{\pi}{k},\]
and $\E^{\P} [|C|\exp(-C)] = \infty$. So~\eqref{hyp:integrability} holds but~\eqref{hyp:finite_entropy} does not. 
\end{itemize}
\end{example}

We have the following observation.
\begin{lemma}
\label{lem:P0_nonempty}
If $\P(C < + \infty) > 0$ then $\mathcal P_0$ is non-empty.
\end{lemma}

\begin{proof}
Under the assumption, there exists an $M > 0$ such that $\P(C < M) > 0$. Let $Z = \1_{\{C \leq M\}} / \P(C \leq M)$. Then $\E^{\P} Z = 1$, so that $d \Q/d \P = Z$ defines a valid probability measure. Furthermore
$\E^{\Q} (C)^+ \leq M < \infty$ and $\E^{\Q} \left[ |\log(Z)| \right] < \infty$. We may conlude that $\Q \in \mathcal P_0$.
\end{proof}
%
%
%

If~\eqref{hyp:integrability} holds, then 
\begin{equation} \label{eq:solution_RE_optimization}
\frac{d\Q^{\star}}{d\P} = Z^{\star} := \frac {\exp(-C)}{\E^{\P} [\exp(-C)]}.
\end{equation} 
defines a probability measure $\Q^{\star}$ that is absolutely continuous with respect to $\P$.

\begin{lemma}
\label{lem:RE_optimization}
Suppose~\eqref{hyp:integrability} holds, and let $\Q^{\star}$ and $Z^{\star}$ be given by~\eqref{eq:solution_RE_optimization}. 
\begin{itemize}
 \item[(i)] $\Q^{\star} \in \mathcal P_0$ if and only if~\eqref{hyp:finite_entropy} holds.
 \item[(ii)] 
 For any $\Q \in \mathcal P_0$, we have that $\Q \ll \Q^{\star}$, and 
\begin{equation} \label{eq:cost_of_other_Q} J(\Q) = \mathcal H(\Q; \Q^{\star}) - \log \E^{\P} \exp(-C), \end{equation}
(so that in particular $J(\Q) = \infty$ if $\mathcal H(\Q;\Q^{\star}) = \infty$);
\item[(ii)] $J$ is strictly convex over $\{ \Q \in \mathcal P_0: \mathcal H(\Q;\Q^{\star}) < \infty\}$ and has unique minimizer $\Q^{\star}$, provided $\Q^{\star} \in \mathcal P_0$.
\end{itemize}
\end{lemma}

\begin{proof}
\begin{itemize}
\item[(i)]
Since $\E^{\P}[ \exp(-C) C \1_{\{ C \geq 0\}}] < \infty$ for any random variable $C$, it remains to check that $\E^{\P}[ \exp(-C) (-C) \1_{\{C < 0\}}] < \infty$, which is the case if and only if~\eqref{hyp:finite_entropy} holds.
\item[(ii)] 
If $\Q$ is not absolutely continuous with respect to $\Q^{\star}$, then there exists a set $E$ of $\Q^{\star}$-measure zero, for which $\Q(E) > 0$. We have $Z^{\star} = 0$, $\P$-almost surely on $E$, so that $C = +\infty$, $\P$-almost surely, on $E$. Therefore $\E^{\Q} (C)^+ \geq \E^{\Q} C \1_E = \infty$. So $\Q \notin \mathcal P_0$.

Now let $\Q \in \mathcal P_0$. We have just seen that $\Q \ll \Q^{\star}$, say with density $Y = \frac{d \Q}{d \Q^{\star}}$. We may choose a version of $Y$ such that $Y = 1$ on $\{ Z^{\star} = 0 \}$.
Then $\Q \ll \P$ with density $Z = \frac{d \Q}{d \P} = Y Z^{\star}$.
Write $K = \E^{\P} \exp(-C)$.
Then
\begin{align*} \mathcal H(\Q;\Q^{\star}) & = \int_{\Omega} Z \log Y \ d \P  = \int_{\Omega \cap \{ Z^{\star} > 0 \} } Z \log Y  \ d \P \\
& = \int_{\Omega  \cap \{ Z^{\star} > 0 \} } Z \left( \log Z - \log Z^{\star} \right) \ d \P = \mathcal H(\Q; \P) - \int_{\Omega  \cap \{ Z^{\star} > 0 \} } Z \log Z^{\star} \ d \P \\
& = \mathcal H(\Q; \P) + \int_{\Omega  \cap \{ Z^{\star} > 0 \} } Z \left( \log K + C \right) \ d \P = \mathcal H(\Q; \P) + \log K + \E^{\P} [Z C].
\end{align*}
In particular
\[ J(\Q) = \E^{\P}  Z C +  \mathcal H(\Q; \P) =  \mathcal H(\Q; \Q^{\star}) -  \log K.\] 
\item[(iii)] Using~\eqref{eq:cost_of_other_Q}, this follows from convexity of relative entropy and the fact that $\mathcal H(\Q; \Q^{\star}) = 0$ if and only if $\Q = \Q^{\star}$.
\end{itemize}
\end{proof}

The following proposition gives a sufficient condition for the existence of a minimizer.

\begin{proposition}
\label{prop:RE_optimization}
Suppose Hypothesis~\eqref{hyp:finite_entropy} holds.
Then $\Q^{\star}$, given by~\eqref{eq:solution_RE_optimization}, is the unique minimizer for  Problem~\ref{prob:REoptimization}, and $J(\Q^{\star}) = - \log \E^{\P} \exp(-C)$.
\end{proposition}

In particular, under~\eqref{hyp:finite_entropy}, $\mathcal H(\Q^{\star};\P) < \infty$ which explains the name of this condition.

\begin{proof}
Note that $\Q^{\star} \in \mathcal P_0$ if~\eqref{hyp:finite_entropy} holds. The stated result is now an immediate consequence of Lemma~\ref{lem:RE_optimization}.
\end{proof}

%
%

The following proposition holds with minimal restrictions on $\P$ and $C$. Since we are primarily interested in the case where Hypothesis~\eqref{hyp:finite_entropy} holds (so as to have the existence result of Proposition~\ref{prop:RE_optimization}), the proof is delegated to the appendix. 
\begin{proposition}
\label{prop:general_optimization}
Suppose $\P(C < + \infty) > 0$. Then
\begin{equation}
\label{eq:general_optimization}
\inf_{\Q \in \mathcal P_0} J(\Q) = - \log \E^{\P} \exp(-C).
\end{equation}
\end{proposition}

To conclude this section, we summarize by distinguishing the following cases, as immediate consequences of Lemma~\ref{lem:RE_optimization}, Proposition~\ref{prop:RE_optimization} and Proposition~\ref{prop:general_optimization}.
\begin{corollary}
\begin{itemize}
 \item[(i)] If Hypothesis~\eqref{hyp:finite_entropy} holds, then there exists a minimizer for Problem~\ref{prob:REoptimization};
 \item[(ii)] If Hypothesis~\eqref{hyp:integrability} holds, but Hypothesis~\eqref{hyp:finite_entropy} does not hold, then Problem~\ref{prob:REoptimization} has optimal value $- \infty < J^{\star} < \infty$, but there does not exist a probability measure $\Q \in \mathcal P_0$ at which this value is attained;
 \item[(iii)] If $\P(C < +\infty) > 0$ but Hypothesis~\eqref{hyp:integrability} does not hold, then Problem~\ref{prob:REoptimization} has optimal value $J^{\star} = - \infty$.
\end{itemize}
\end{corollary}

\subsection{Notes and remarks}

\begin{remark}
One may wish to include a factor $\beta > 0$ in the problem formulation, to indicate the relative importance of minimizing $\E^{\Q} C$ compared to minimizing $\mathcal H(\Q ; \P)$ to obtain the form $J(\Q) =  \E^{\Q} C + \frac 1 {\beta} \mathcal H(\Q ; \P)$. In this case the optimal value becomes $J^{\star} = - \frac 1 {\beta} \log \E^{\P} \exp(-\beta C)$, so that $\beta$ admits the interpration of inverse temperature, in the context of statistical physics. 
 \end{remark}
 
 \begin{remark}
 \label{rem:sufficient_for_finite_entropy}
 A sufficient condition for Hypothesis~\eqref{hyp:finite_entropy} to hold is that for some $\gamma > 1$,
 \begin{equation} \label{eq:sufficient_for_finite_entropy} 0 < \E^{\P} \exp(-\gamma C) < \infty.\end{equation}
 Indeed, if this is the case, then, since $x \leq \frac 1 {\varepsilon e} \exp(\varepsilon x)$ for $x \geq 0$,
 \[ \E^{\P} \left[ \exp(-C) (-C) \1_{\{C \leq 0\}}\right] \leq \frac 1 {(\gamma - 1) e} \E^{\P} \left[ \exp(-\gamma C) \1_{\{ C \leq 0 \}} \right] < \infty.\]
 (And, as always, $\E^{\P} \left[ \exp(-C) |C| \1_{\{ C > 0\}} \right] \leq \frac 1 e$.)
 In turn~\eqref{eq:sufficient_for_finite_entropy}, and therefore~\eqref{hyp:finite_entropy}, are implied by the conditions that $\P(C > K) = 1$ for some $K \in \R$, and $\P(C = \infty) < 1$.
 \end{remark}

\section{Relative entropy weighted optimization with respect to Wiener measure}
\label{sec:dynamic_RE_optimization}

Our primary interest lies in the theory discussed in Section~\ref{sec:RE_optimization} applied to the special case where all randomness is generated by a $d$-dimensional Wiener process. 

Equivalent changes of measure (satisfying mild conditions) may in this case be expressed as a Girsanov type transformation. The corresponding change of drift will constitute the `control process'. A crucial observation is that the relative entropy of such a transformation is given by the squared control costs. However, quite importantly, we will not restrict ourselves to the case in which the optimal measure is equivalent to the Wiener measure, and will only require absolute continuity, with finite relative entropy. This allows us to consider cost functionals that may have value $+\infty$ with positive probability.

In Section~\ref{sec:general_results} we will discuss the general setting. In Section~\ref{sec:control_as_malliavin_derivative} we specialize to the case in which the cost functional has a Malliavin derivative, which provides a convenient expression for the optimal drift.


\subsection{General results}
\label{sec:general_results}

\begin{definition}[Canonical Wiener process]
Let $\Omega = C([0,\infty);\R^d)$, i.e. the space of continuous functions mapping $[0,\infty)$ into $\R^d$. For $t \geq 0$, let $\mathcal F_t^o = \sigma \left( \left\{ \omega(s) : 0 \leq s \leq t \right\} \right)$ for $t \geq 0$. Let $\P$ denote Wiener measure, and for $t \geq 0$, let $(\mathcal F_t)$ be the right-continuous completion of $(\mathcal F_t^o)$ with respect to $\P$ (see \cite[Section II.67]{RogersWilliams1994a}). Let $X_t(\omega) := \omega(t)$ for $\omega \in \Omega$ and $t \geq 0$, so that $X$ is a standard Brownian motion in $\R^d$ under the probability measure $\P$. Let $\mathcal F := \mathcal F_{\infty} := \sigma(\cup_{t \geq 0} \mathcal F_t)$. The collection $(\Omega, (\mathcal F_t)_{t \geq 0}, \mathcal F, \P, X)$ will be referred to as a \emph{canonical $d$-dimensional Wiener process}.
\end{definition}

Throughout this section let $(\Omega, (\mathcal F_t)_{t \geq 0}, \mathcal F, \P,X)$ denote a canonical $d$-dimensional Wiener process. Furthermore let $C : \Omega \rightarrow \R$ satisfy~\eqref{hyp:finite_entropy}, i.e. 
\[ \P(C < + \infty) > 0 \quad \mbox{and} \quad \E^{\P} \left[ \exp(-C) |C|\right] < \infty.\]

The following result and its practical application (Remark~\ref{rem:explicit_follmer_drift}) are of fundamental importance in the explicit computation of the optimal drift in general cases.

\begin{proposition}[F\"ollmer's drift]
\label{prop:follmer_drift}
Let $\Q$ be a probability measure on $(\Omega, \mathcal F)$ that is absolutely continuous with respect to $\P$. There exists a adapted process $(U_t)$ such that under $\Q$ the following hold:
\begin{itemize}
 \item[(i)] The process $X^U_t = X_t - \int_0^t U_s \ ds$ is a $\Q$-Brownian motion;
 \item[(ii)] $\mathcal H(\Q;\P) = \half \E^{\Q} \int_0^{\infty} |U_t|^2 \ d t$.
\end{itemize}
\end{proposition}

\begin{proof}
See \cite[Section 2]{Follmer1986} or \cite[Theorem 2]{Lehec2011}.
\end{proof}

\begin{remark} \label{rem:explicit_follmer_drift}
Since our interest is in obtaining explicit expressions for the optimal feedback control, we also state here the construction of the F\"ollmer drift $(U_t)$, and hence of the feedback control $u$.
Let $Z = \frac{d \Q}{d \P}$, and write $Z_t = \E^{\P} [ Z \mid \mathcal F_t]$, $t \geq 0$. Since $(Z_t)$ is a continuous martingale, it allows the martingale representation $Z_t = 1 + \int_0^t \langle V_s, \ d X_s \rangle$ \cite[Theorem IV.36.1]{RogersWilliams1994b}. The F\"ollmer drift is given by $U_t := (V_t / Z_t) \1_{Z_t > 0}$. See \cite[Section 2]{Follmer1986} or \cite[Theorem 2]{Lehec2011} for details. 
\end{remark}

\begin{definition}[Feedback control]
A mapping $u : [0,\infty) \times C([0,\infty);\R^d) \rightarrow \R^d$ is called a \emph{feedback control} if the process $U_t := u(t,X)$ is $(\mathcal F_t)$-adapted.
\end{definition}

Suppose $u$ is a feedback control and consider the stochastic differential equation
\begin{equation} \label{eq:SDE} Y_t = \int_0^t u(s,Y) \ d s + B_t \quad \mbox{for} \quad t \geq 0.\end{equation}
Recall that a \emph{weak solution} of~\eqref{eq:SDE} consists of a filtered probability space $(\widetilde \Omega, \widetilde {\mathcal F}, (\widetilde {\mathcal F}_t), \widetilde \P)$ satisfying the usual conditions, on which is defined a $d$-dimensional standard Brownian motion along with a $(\mathcal F_t)$-adapted process $Y$ with values in $\R^d$ such that~\eqref{eq:SDE} is satisfied. Note that $(Y_t)$ is not necessarily adapted to the filtration generated by $(B_t)$.

The following result shows how changes of measure may be associated with feedback controls and vice versa.

\begin{proposition}
\label{prop:feedback_control}
\begin{itemize}
 \item[(i)] Let $u$ be a feedback control and suppose $(\widetilde \Omega, \widetilde {\mathcal F}, ( \widetilde {\mathcal F}_t), \widetilde \P, Y, B)$ is a weak solution to~\eqref{eq:SDE}, such that $\E^{\widetilde \P} \left[ \int_0^{\infty} |u(t, Y)|^2 \ d t \right] < \infty$. Then there exists a probability measure $\Q \ll \P$ on $(\Omega, \mathcal F)$, such that under $\Q$, the coordinate process has the same distribution as $Y$, and such that $\mathcal H(\Q;\P) = \frac 1 2 \E^{\widetilde \P} \int_0^{\infty} |u(t,Y)|^2 \ d t$.
 \item[(ii)] Suppose $\Q \ll \P$ with $\mathcal H(\Q;\P) < \infty$. Then there exists a feedback control $u$ and adapted processes $Y$ and $B$ on $(\Omega, \mathcal F, (\mathcal F_t))$ satisfying~\eqref{eq:SDE} such that $\Q$ is the law of $Y$, $B$ is a standard Brownian motion under $\Q$, and $\mathcal H(\Q;\P) = \half \E^{\P} \int_0^{\infty} |u(t,Y)|^2 \ d t$.
\end{itemize}
\end{proposition}

\begin{proof}
\begin{itemize}
 \item[(i)] See e.g. \cite[Proposition 1]{Lehec2011}. 
 \item[(ii)] Let $U$ be the F\"ollmer drift corresponding to $\Q$, and note that $u(t,X) := U_t(\omega)$ is a feedback control. Furthermore $B_t := X_t - \int_0^t u(s,X) \ d s$ is a standard Brownian motion under $\Q$, so that $Y := X$ satisfies~\eqref{eq:SDE}.
\end{itemize}
In either case Proposition~\ref{prop:follmer_drift}(ii) gives the stated expression for the entropy.
\end{proof}

%
%
%
%
 
The following result is now a simple combination of Proposition~\ref{prop:RE_optimization} and Proposition~\ref{prop:feedback_control}(ii).

\begin{corollary}
\label{cor:dynamic_RE_optimization}
Suppose $(\Omega, \mathcal F, (\mathcal F_t), \P, X)$ is a $d$-dimensional canonical Wiener process.
Let $C : \Omega \rightarrow \R$ be $\mathcal F$-measurable such that~\eqref{hyp:finite_entropy} holds.
Then
there exists a probability measure $\Q^{\star} \in \mathcal P_0$ that solves Problem~\ref{prob:REoptimization}.
Furthermore there exists a feedback control $u$, such that 
\begin{itemize}
\item[(i)] $U_t = u(t,X)$ is the F\"ollmer drift (of Proposition~\ref{prop:follmer_drift}) associated with $\Q^{\star}$; 
\item[(ii)] Under $\Q^{\star}$, the law of the coordinate process $X$ coincides with the law of the solution $Y$ of~\eqref{eq:SDE}; and
\item[(iii)] $\half \E^{\Q^{\star}} \int_0^{\infty} |u(s,X)|^2 \ ds  = \mathcal H(\Q^{\star};\P) < \infty$.
\end{itemize}
\end{corollary}

%
\begin{remark}
The practical interpretation of Corollary~\ref{cor:dynamic_RE_optimization} is that we solve the control problem
\[ \inf_u \E^{\widetilde \P} \left[ C(Y^u) + \half \int_0^{\infty} |u(s,Y^u)|^2 \ d s \right],\]
with the infimum taken over all feedback controls $u$, and where $Y^u$ satisfies~\eqref{eq:SDE} with respect to a $\widetilde \P$-Brownian motion $B$ on some fixed probability space $(\widetilde \Omega, \widetilde \P)$. The problem with making this interpretation rigourous, is that one would have to impose conditions on $u$ such that a strong solution to~\eqref{eq:SDE} always exists. As illustrated by Corollary~\ref{cor:dynamic_RE_optimization}, we are able to obtain, in general, a weak solution to the mentioned problem. In \cite[Section 2.3]{Lehec2011} conditions are discussed such that a strong solution exists. 
\end{remark}

\begin{remark}
The result of Corollary~\ref{cor:dynamic_RE_optimization} is relatively straightforward if $\P(C = +\infty) = 0$. In this case the optimal measure $\Q^{\star}$ with density $Z^{\star} = \frac{d \Q^{\star}}{d \P} = \frac{\exp(-C)}{\E^{\P} \exp(-C)}$ is equivalent to $\P$, so that standard Girsanov theory (e.g. \cite[Theorem IV.38.5(i)]{RogersWilliams1994b}) can be applied to find the representation
\begin{equation}
\label{eq:girsanov_density} Z_t^{\star} =\exp\left( \int_0^t \langle u(s,X), \ dX_s \rangle - \half \int_0^t |u(s,X)|^2 \ d s \right)
\end{equation} for the density process $Z_t^{\star} := \E^{\P}[Z^{\star} |\mathcal F_t]$;
it only remains to verify the stated expression for the relative entropy, which boils down to a simple computation.
We have to be more careful if $\P(C = +\infty) > 0$. In this case $\P(Z^{\star} = 0) > 0$ so that the measures $\Q^{\star}$ and $\P$ are not equivalent. The representation~\eqref{eq:girsanov_density} is still valid but only $\Q^{\star}$-almost surely; see \cite[Section 2]{Follmer1985} or \cite[Theorem 2]{Lehec2011}.
\end{remark}

The following proposition provides a method to compute the optimal drift through Monte Carlo sampling, by employing the law of large numbers to approximate the (conditional) expectation. In Section~\ref{sec:minimize_maximum_monte_carlo} this idea is used in the approximate minimization of the running maximum of a Brownian motion with drift.

\begin{proposition}
\label{prop:control_as_derivative}
Let $\Q \ll \P$ with $\mathcal H(\Q;\P) < \infty$ with density $Z = \frac{d \Q}{d \P}$ and let $u$ be the feedback control of Proposition~\ref{prop:feedback_control}. Suppose that $Z$ is $\mathcal F_T$-measurable for some $T \geq 0$.
Then the corresponding F\"ollmer drift $U_t = u(t,X)$ at time $t$ satisfies
\begin{equation}
\label{eq:control_as_derivative}
U_t =  \lim_{h \downarrow 0} \frac 1 {h} \E^{\Q}[ (X_{t+h}-X_t) \mid \mathcal F_t] \quad \mbox{in $L^2(\Q)$ for $0 \leq t \leq T$.}
\end{equation}

\end{proposition}

\begin{proof}
This is \cite[Proposition 2.5]{Follmer1986}.
\end{proof}

\subsubsection{Examples and remarks}
\label{sec:examples}
In all the examples below, $Z^{\star} = \frac{d \Q^{\star}}{d \P}$, with $\Q^{\star}$ as in Corollary~\ref{cor:dynamic_RE_optimization}. Also, $Z_t^{\star} = \E^{\P} \left[ Z^{\star} \mid \mathcal F_t \right]$, $t \geq 0$.
\begin{remark}
Finite time horizon problems fit well within the theory of this section: In this case one just considers a cost random variable $C$ that is $\mathcal F_T$-measurable for some deterministic $T > 0$. The optimal drift $(U_t)$ will be $\mathcal F_T$-measurable for all $t \geq 0$, by its construction. 
\end{remark}
As an illustration of this remark, consider the following example.


\begin{example}[Constrained problem, obstacle avoidance]
 Let $A \subset \R^d$ be open and let $T > 0$. Then the event
 $\{ \omega: X_T \in A\}$ is $\mathcal F_T$-measurable.
 Let 
 \[ C = \left\{ \begin{array}{ll}+ \infty \quad & \mbox{if $X_T \in A$}, \\
 0 \quad & \mbox{otherwise.} \end{array} \right.\]
 Then~\eqref{hyp:finite_entropy} is satisfied as long as $\P(C < + \infty) = \P(X_T \notin A) > 0$. 
 Optimally controlled paths will avoid the set $A$ at time $T$.
For example, in the one-dimensional case, let us take $A = (a,b)$, with $a < b$, possibly with $a = -\infty$ or $b = + \infty$, but not both.
We compute for $t < T$,
\[ N_t  := \E^{\P} \left[\exp(-C) \mid \mathcal F_t \right] = \E^{\P} \left[ \1_{X_T \notin A} \mid \mathcal F_t \right] 
 = 1 - \frac 1 {\sqrt{2 \pi (T-t)}}  \int_a^b  \exp\left( -(x-X_t) ^2/2(T-t) \right)  \ d x.\]
By It\^o calculus, differentiating under the integral, for $t < T$,
\begin{align*}
 dN_t & = - \left( \frac 1 {\sqrt{2 \pi (T-t)}} \int_a^b \left\{ \exp\left( -(x-X_t) ^2/2(T-t) \right) \left( \frac{x - X_t}{T-t} \right) \right\} \ d x \right) \ d X_t  \\
 &  = \left( \frac {1} {\sqrt{2 \pi (T-t)}} \left\{  \exp \left( -(b - X_t)^2/2 (T-t) \right) - \exp \left( -(a - X_t)^2/2 (T-t) \right) \right\} \right) \ d X_t
\end{align*}
Since $Z_t^{\star} = N_t / N_0$, and by Remark~\ref{rem:explicit_follmer_drift}, we find that
\[ U_t = \frac{ \exp \left( -(b - X_t)^2/2 (T-t) \right) - \exp \left( -(a - X_t)^2/2 (T-t) \right)}{ \sqrt{2 \pi (T-t)} - \int_a^b  \exp\left( -(x-X_t) ^2/2(T-t) \right)  \ d x } \1_{t < T} \]
is the optimal drift $(U_t)$.
 \end{example}

Note that we have no difficulty in dealing with the fact that $\P(C = +\infty) > 0$. This is also illustrated by the following example.

\begin{example}[Whittle's `flypaper' example {\cite[Section 20.3]{Whittle1982},\cite[Sections V.7, V.15]{RogersWilliams1994b}}]
\label{ex:fly_paper}
Consider the case in which $X$ is a one-dimensional Wiener process, i.e. $d = 1$. Let $\tau = \inf \{ t  > 0 : X_t \in \{ a, b \} \}$, with $a < 0 < b$. As is well known, $\E \tau < \infty$. Consider the cost functional $C = k_a \1_{X_{\tau} = a} + k_b \1_{X_{\tau} = b} + \gamma \tau$, with $k_a$ and $k_b$ given scalar constants in $(-\infty, \infty]$, not both equal to $+\infty$, and $\gamma \geq 0$. Then~\eqref{hyp:finite_entropy} is satisfied for $C$.
Suppose that $\varphi$ is the unique twice continuously differentiable solution on $[a,b]$ to the linear differential equation
\[ \left\{ \begin{array}{ll} \half \varphi''(x) - \gamma \varphi(x) = 0, \quad & (a < x < b), \\
            \varphi(a) = \exp(-k_a),  \quad &  \varphi(b) = \exp(-k_b).
           \end{array} \right.
\]
Define $N_t := \exp(-\gamma(t \wedge \tau)) \varphi(X_{t \wedge \tau})$. Then by It\^o's formula, for $0 \leq t \leq \tau$,
\[ dN_t = - \gamma N_t \ d t + \half \exp(-\gamma t) \varphi''(X_t) + \exp(-\gamma t) \varphi'(X_t) \ d X_t = \exp(-\gamma t) \varphi'(X_t) \ d X_t,\]
and $d N_t = 0$ for $t > \tau$,
so that $(N_t)$ is a local martingale. 
In fact, since $\varphi$ has bounded derivative, $(N_t)$ is a martingale. Since $N_{\infty} = \exp(-C)$, we have $Z_t^{\star} = N_t / N_0$. By Remark~\ref{rem:explicit_follmer_drift}, it follows that 
\[ u(t,X) = \varphi'(X_t) / \varphi(X_t) \1_{t < \tau} \]
 is the optimal feedback control.
%
\end{example}

\begin{remark}
\label{rem:control_of_diffusion}
Example~\ref{ex:fly_paper} illustrates the basic approach to more general problems where $C$ is of the form 
\[ C = \int_0^{\tau} v(s, \xi_s) \ d s + k(\tau, \xi_\tau),\]
with $(\xi_t)$ a process in $\R^n$ satisfying an SDE of the form
\begin{equation} \label{eq:Ito-diffusion} \xi_t = \xi_0 + \int_0^t b(s,\xi_s) \ d s + \int_0^t \sigma(s,\xi_s) \ d X_s,\end{equation}
and $\tau = \inf \{ t > 0 : \xi_t \notin U \ \mbox{or} \ t > T \}$ denoting the exit time of a region $U \times [0,T)$ with $U$ open and $T \in (0, \infty]$. Assume that $\P(\tau < \infty) = 1$. The optimal solution of Corollary~\ref{cor:dynamic_RE_optimization} may in principle be obtained by finding a solution to the PDE
\[ \label{eq:linear_pde} \left\{ \begin{array}{ll} L \varphi(t,x) - v(t,x) \varphi(t,y) = 0, \quad & x \in U, 0 < t < T,  \\
            \varphi(t,x) = \exp(-k(t,x)), \quad & x \in \partial U \ \mbox{or} \ t = T,
           \end{array} \right.
\]
where $L$ is the infinitesemal generator of~\eqref{eq:Ito-diffusion},
and verifying that $N_t := \exp \left( - \int_0^{t \wedge \tau} v(s,\xi_s) \ d s \right) \varphi(t \wedge \tau,\xi_{t \wedge \tau})$
is a martingale. Then $Z_t^{\star} = N_t/N_0$, and an application of It\^o's formula provides (by Remark~\ref{rem:explicit_follmer_drift}) the optimal drift 
\[ U_t = \frac {  [\sigma(t, \xi_t)]^T \nabla \varphi(t,\xi_t)}{\varphi(t, \xi_t)} \1_{t \leq \tau} \1_{\varphi(t,\xi_t) > 0 }.\]
The optimal drift is then a Markov control for $(\xi_t)$ (i.e. of the form $u(t, \xi_t)$), and under the optimal change of measure, $\xi$ satisfies (in law) the `controlled' SDE
\[ \xi_t = \xi_0 + \int_0^t \left\{ b(s,\xi_s) + \sigma(s, \xi_s) u(s, \xi_s)\right\}  \ d s + \int_0^t \sigma(s,\xi_s) \ d B_s,\]
where $(B_t)$ is a standard Brownian motion.
\end{remark}

\begin{remark}
In \cite{BoueDupuis1998}, it is proven that for $C$ bounded from above, and $\mathcal F_1$-measurable (i.e. in the finite time horizon case), we have
\[ \inf_{U} J(U) = -\log \E^{\P} \exp(-C),\]
where the minimization is over a suitable class of processes that are progressively measurable with respect to $(\mathcal F_t)$.
This result is used in the field of large deviations theory \cite{BoueDupuis1998, DupuisEllis1997}.  See also \cite[Theorems 4, 9]{Lehec2011} for an elegant alternative approach, extending the above Bou\'e-Dupuis formula to infinite time horizon.
Our focus is on explicit computations in optimal control, so we have taken  care so that our approach
%
(i) is not restricted to finite time horizon problems; (ii) establishes explicitly the existence of a minimizer; and (iii) does not impose the restriction on the cost functional to be bounded from above.
\end{remark}

%

\begin{remark}
Our set up, using Wiener measure on the time interval $[0,\infty)$ allows for control over an infinite time horizon. At the same time, by Corollary~\ref{cor:dynamic_RE_optimization},  condition~\eqref{hyp:finite_entropy} excludes e.g. the case in which the control process is $\Q$-almost surely equal to a non-zero constant on an infinite time horizon (since then $\mathcal H(\Q;\P) = \infty$).
However, condition~\eqref{hyp:finite_entropy} allows us to consider problems whose cost function is $\mathcal F_{\tau}$-measurable, with $\tau$ a stopping time that is unbounded but $\P$-almost surely finite, as illustrated by Example~\ref{ex:fly_paper}.
\end{remark}

\begin{remark}
An alternative method for finding an expression for the optimal drift in the case of the control of diffusion processes with cost function as in Remark~\ref{rem:control_of_diffusion} is through the dynamic programming principle, i.e. the Hamilton-Jacobi-Bellman PDE. For the class of problems considered in this paper, the HJB equation may be transformed into the linear PDE~\eqref{eq:linear_pde} through the \emph{Hopf-} or \emph{logarithmic transform} \cite{FlemingSoner2009, Kappen2005}.
\end{remark}

\subsection{Computation of optimal drift by means of Clark-Ocone representation}
\label{sec:control_as_malliavin_derivative}
In this section we will apply the Clark-Ocone theorem of Malliavin calculus to obtain an explicit representation of the optimal drift in terms of the Malliavin derivative of the cost random variable. This representation also occurs e.g. in \cite[Lemma 3]{Lehec2011}.
Section~\ref{sec:minimize_maximum} provides an illustration of the use of Corollary~\ref{cor:control_as_malliavin_derivative}.

As usual, let $(\Omega, (\mathcal F_t)_{t \geq 0}, \mathcal F,\P, X)$ denote a canonical $d$-dimensional Wiener process.
Recall the following definitions and notations (see \cite[Chapter 1]{Nualart2006} for details).
Let $H$ denote the Hilbert space $L^2([0,\infty); \R^d)$. For suitable random variables $F : \Omega \rightarrow \R$, the Malliavin derivative $D F \in L^2((\Omega, \P); H)$ may be defined. The domain of $D$ in $L^2(\Omega, \mathcal F, \P)$ is denoted by $\mathbb D^2$.
%
%
%
If $F \in \mathbb D^2$, so that $D F \in L^2(\Omega;H) \cong L^2(\Omega \times [0,\infty);\R^d)$, then we may identify $DF$ with a stochastic process, denoted $D_t F$, as usual:
\[ (D_t F)(\omega) = DF(\omega,t), \quad \omega \in \Omega, \quad t \geq 0.\]
Note that $(D_t F)_{t \geq 0}$ is not necessarily adapted, and the value of $D_t F(\omega)$ is defined $\P \otimes \mathrm{Leb}$-almost everywhere.
We can now quote the following result, which is essentially an application of the Clark-Ocone formula \cite[Proposition 1.3.14]{Nualart2006}.

\begin{lemma}
\label{lem:control_as_malliavin_derivative}
Suppose $Z \in \mathbb D^2$ and let $(U_t)$ be the F\"ollmer process associated to the measure $\Q$ with density $Z$ with respect to $\P$. Let $Z_t := \E^{\P} \left[ Z \mid \mathcal F_t \right]$ for $t \geq 0$. Then $\P \otimes \mathrm{Leb}$-almost everywhere,
\[ U_t = \frac{\E^{\P} \left[ D_t Z \mid \mathcal F_t \right]}{Z_t} \1_{\{ Z_t > 0 \}}.\]
\end{lemma}

\begin{proof}
See \cite[Lemma 3]{Lehec2011}.
\end{proof}

\begin{corollary}
\label{cor:control_as_malliavin_derivative}
Suppose $\P(0 \leq C < \infty) = 1$, and $C \in \mathbb D^{2}$. 
Then the optimal drift $(U_t)_{t \geq 0}$ of Corollary~\ref{cor:dynamic_RE_optimization} satisfies, $\P \otimes \mathrm{Leb}$-almost everywhere,
\[ U_t = - \frac{\E^{\P} \left[ \exp(-C) D_t C \mid \mathcal F_t \right]}{\E^{\P} \left[ \exp(- C) \mid \mathcal F_t \right] }.\]
\end{corollary}

\begin{proof}
Hypothesis~\eqref{hyp:finite_entropy} is satisfied. By Proposition~\ref{prop:RE_optimization}, the optimal density in Problem~\ref{prob:REoptimization} has density $Z^{\star} = \frac{\exp(-C)}{\E^{\P}[\exp(-C)]}$. Since $\P(C < +\infty) = 1$, $Z > 0$, $\P$-almost surely and therefore $Z_t > 0$, $\P$-almost surely for all $t \geq 0$ \cite[Theorem II.78.1(ii)]{RogersWilliams1994a}. Furthermore $x \mapsto \exp(-x)$ has bounded derivative for $x$ bounded from below. Therefore we may apply the chain rule of Malliavin calculus to $\exp(-C)$ (see \cite[Proposition 1.2.3]{Nualart2006}), to conclude that $\exp(-C) \in \mathbb D^{2}$ with $D_t \exp(-C) = - \exp(-C) D_t C$. The result now follows from Lemma~\ref{lem:control_as_malliavin_derivative}. 
\end{proof}

%

%

\begin{remark}
In \cite{OconeKaratzas1991}, the Clark-Ocone formula is used in a different optimization context, namely for portfolio optimization in mathematical finance.
\end{remark}

\subsection{Minimization of the maximum of a Brownian motion with drift over a finite time horizon}
\label{sec:minimize_maximum}
In this section we illustrate the theory by obtaining the optimal drift for the minimization of the maximum of a standard Brownian motion with drift over a finite time horizon. As it turns out, the value of optimal drift can be explicitly computed as a function of the difference between the running maximum and the current value of the Wiener process. 
For a different optimization problem related to the maximum of a Wiener process see \cite{HeinricherStockbridge1991}.

First, in Section~\ref{sec:minimize_maximum_malliavin}, we apply the result of Section~\ref{sec:control_as_malliavin_derivative} to obtain an explicit expression for the optimal drift through the Clark-Ocone formula. Next, in Section~\ref{sec:minimize_maximum_monte_carlo}, we show how to apply Proposition~\ref{prop:control_as_derivative} to obtain a Monte-Carlo estimate of the optimal drift.

Let $(\Omega, (\mathcal F_t)_{t \geq 0}, \mathcal F, \P, X)$ denote a one-dimensional canonical Wiener process. Define $M_t := \max_{0 \leq s \leq t} X_s$ and take $C := M_T$ for some $T > 0$.
Hypothesis~\eqref{hyp:finite_entropy} is satisfied by Remark~\ref{rem:sufficient_for_finite_entropy} and the observation that $\P(M_T = \infty) = 0$.

\subsubsection{Exact solution}
\label{sec:minimize_maximum_malliavin}
For the distribution of $M_t$ we have by virtue of the reflection principle \cite[Section 2.8.A]{KaratzasShreve1991}
\[ \P \left( M_t \geq a \right) = \left(\frac 2 \pi \right)^{1/2} \int_{a / \sqrt{t}}^{\infty} \exp(-\xi^2/2) \ d \xi, \quad t \geq 0, a \geq 0.\]
We will make use of the error function (erf) and complimentary error function (erfc), defined by
\[ \erf(x) := \frac 2 {\sqrt{\pi}} \int_0^x \exp(-\eta^2) \ d \eta, \quad \erfc(x) := 1 - \erf(x) = \frac 2 {\sqrt{\pi}} \int_x^{\infty} \exp(-\eta^2) \ d \eta, \quad x \geq 0.\]

We compute the optimal drift corresponding to the minimization of $C$, using Corollary~\ref{cor:control_as_malliavin_derivative}. This means we need to compute
$\E^{\P} \left[ \exp(-M_T) D_t M_T \mid \mathcal F_t \right]$ and $\E^{\P} \left[ \exp(-M_T) \mid \mathcal F_t \right]$.
We start with the latter. Conditional on $\mathcal F_t$, the event $M_T = M_t$ occurs when the maximum over $[t,T]$ does not exceed $y := M_t$. This has the same probability as the event that the maximum over $[0,T-t]$ does not exceed $y - X_t$, so
\[ \P(M_T = M_t | X_t = x, M_t = y) = \P(M_{T-t} \leq y - x) = \left(\frac 2 \pi \right)^{1/2}\int_0^{\frac{M_t -X_t}{\sqrt{T-t}}} \exp( - \xi^2/2) \ d \xi.\]
For $0 \leq x \leq y < z$ we compute
\begin{align*}
\P \left( M_T \geq z | X_t = x, M_t = y\right) &= \P \left(M_T = M_t | X_t = x, M_t = y\right) + \P \left( M_{T-t} \geq z - x \right) \\
& =\left(\frac 2 \pi \right)^{1/2}\int_0^{\frac{M_t-X_t}{\sqrt{T-t}}} \exp( - \xi^2/2) \ d \xi + \left(\frac 2 \pi \right)^{1/2} \int_{\frac{z-X_t}{\sqrt{T-t}}}^{\infty} \exp(-\xi^2/2) \ d \xi.
\end{align*}
Therefore the density function of $M_T$ conditional on $\mathcal F_t$ is equal to
\[ f_{M_T | \mathcal F_t}(\xi) = \left(\frac{2}{\pi(T-t)}\right)^{1/2} \exp\left( - \frac{(\xi-X_t)^2}{2 (T-t)} \right), \quad \mbox{for} \ \xi > M_t \geq X_t.\]
We compute
\begin{align}
\nonumber 
N_t := & \E^{\P} \left[ \exp(- M_T) \mid \mathcal F_t \right] = \E^{\P} \left[ \exp(- M_T) \1_{M_T = M_t} \mid \mathcal F_t \right] + \E^{\P} \left[ \exp(- M_T) \1_{M_T > M_t}\mid\mathcal F_t \right] \\
\nonumber & = \exp(- M_t) \P(M_T = M_t\mid\mathcal F_t) + \E^{\P} \left[ \exp(- M_T) \1_{M_T > M_t}\mid\mathcal F_t \right] \\
\nonumber & = \exp(- M_t) \left(\frac 2 \pi \right)^{1/2}\int_0^{\frac{M_t-X_t}{\sqrt{T-t}}} \exp( - \xi^2/2) \ d \xi \\
\nonumber & \quad \quad + \left( \frac 2 {\pi(T-t)} \right)^{1/2} \int_{M_t}^{\infty} \exp(- \xi) \exp\left( - \frac{(\xi-X_t)^2}{2 (T-t)} \right) \ d \xi \\
\label{eq:N_t_maximum}& = \exp(- M_t) \erf \left( \frac{M_t - X_t}{\sqrt{2 (T-t)}} \right) + \exp \left( - X_t + \half (T-t) \right) \erfc \left( \frac{M_t - X_t + (T-t)}{\sqrt{2 (T-t)}} \right)
\end{align}
The Malliavin derivative of $M_T$ is given by $D_t M_T = \1_{[0,\tau]}(t) = \1_{M_t < M_T}$, where $\tau$ is the a.s. unique point where $X$ attains its maximum \cite[Exercise 1.2.11]{Nualart2006}. 
%
Therefore
\begin{align*} V_t & := - \E^{\P} \left[ \exp(-M_T ) D_t M_T \mid \mathcal F_t \right]  =- \E^{\P} \left[ \exp(- M_T) \1_{ M_T > M_t } \mid \mathcal F_t \right] \\
 & = -\left( \frac{2}{\pi(T-t)}\right)^{1/2} \int_{M_t}^{\infty}  \exp \left( - \xi - \frac{(\xi - X_t)^2}{2 (T-t)} \right) \ d \xi \\
& = - \exp \left( -  X_t + \half  (T-t) \right) \erfc \left( \frac{M_t - X_t + (T-t)}{\left( 2 (T-t) \right)^{1/2} } \right).
\end{align*}
We finally compute the optimal feedback control 
\begin{equation} \label{eq:solution_1} u(t, X) := V_t / N_t = u(t, X_t, M_t), \end{equation} where (with some abuse of notation), for $0 \leq t < T$, $x \in \R$, and $m \geq x$,
\begin{equation} \label{eq:solution_2} u(t, x, m) = \frac {-\exp \left( - x + \half (T-t) \right) \erfc \left( \frac{m - x + T-t}{\sqrt{2 (T-t)}} \right)}{\exp(- m) \erf \left(\frac {m - x}{\sqrt{ 2 (T-t)}}\right) 
+ \exp\left( -  x + \half (T-t) \right) \erfc \left( \frac {m - x + T-t}{\sqrt{2 (T-t)}} \right)     
}.\end{equation}
\begin{remark}
This example illustrates how the theory developed in this paper may be applied to obtain non-Markov controls, and therefore provides a method that applies where a dynamic programming (i.e. the HJB equation) can not be used. One could in principle obtain a Markov control for the extended process $(X_t, M_t)$. However solution by means of dynamic programming is still far from straightforward.
\end{remark}

\begin{remark}
Alternatively, one can find the required martingale representation by applying It\^o's formula to the expression for $N_t$ obtained in~\eqref{eq:N_t_maximum}, as in the examples of Section~\ref{sec:examples}. Since the semimartingale decomposition of $Z_t^{\star} := N_t/N_0$ is unique \cite[Section IV.31]{RogersWilliams1994b}, the decomposition obtained using It\^o's formula is necessarily equal to the required martingale representation of $(Z_t^{\star})$ as in Remark~\ref{rem:explicit_follmer_drift}.
\end{remark}
%

\subsubsection{Approximate solution through Monte-Carlo sampling}
\label{sec:minimize_maximum_monte_carlo}
\begin{figure}
\begin{center}
\begin{subfigure}[b]{0.45 \textwidth}\includegraphics[width=\textwidth]{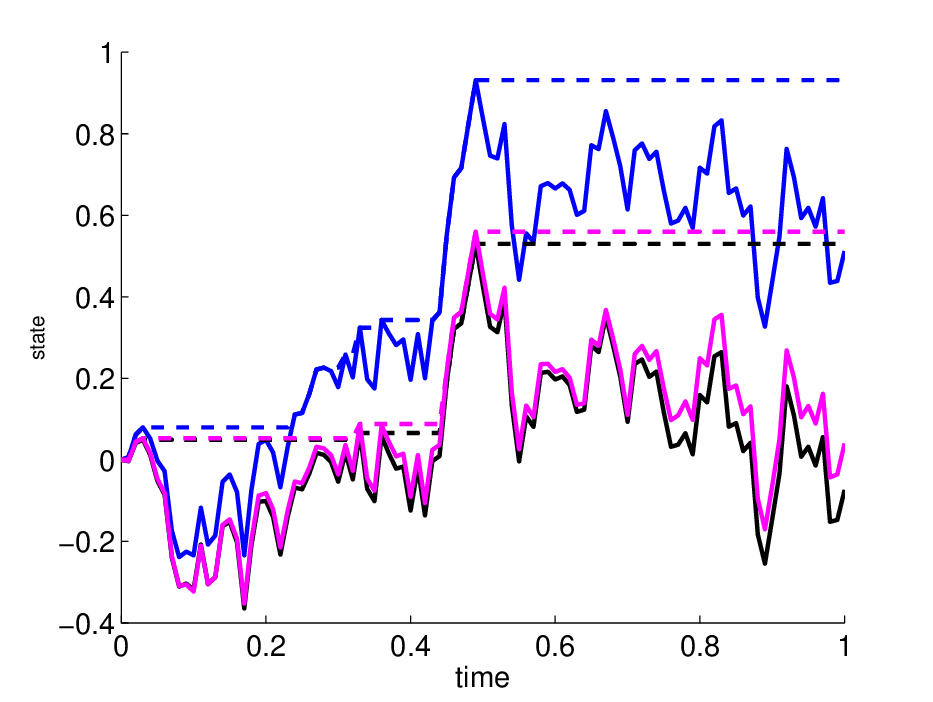}
\caption{Sample path of Brownian motion with drift (solid line) along with its maximum process (dashed line), uncontrolled (blue), controlled with exact optimal drift (black) and controlled with approximately optimal drift, obtained by Algorithm~\ref{alg:monte-carlo} (magenta). It is seen how the controlled paths reduce the value of the maximum process at time $t = 1$, where the exact control outperforms the Monte-Carlo solution (but possibly with higher control cost).}
\end{subfigure}
\hfill
\begin{subfigure}[b]{0.45 \textwidth}\includegraphics[width=\textwidth]{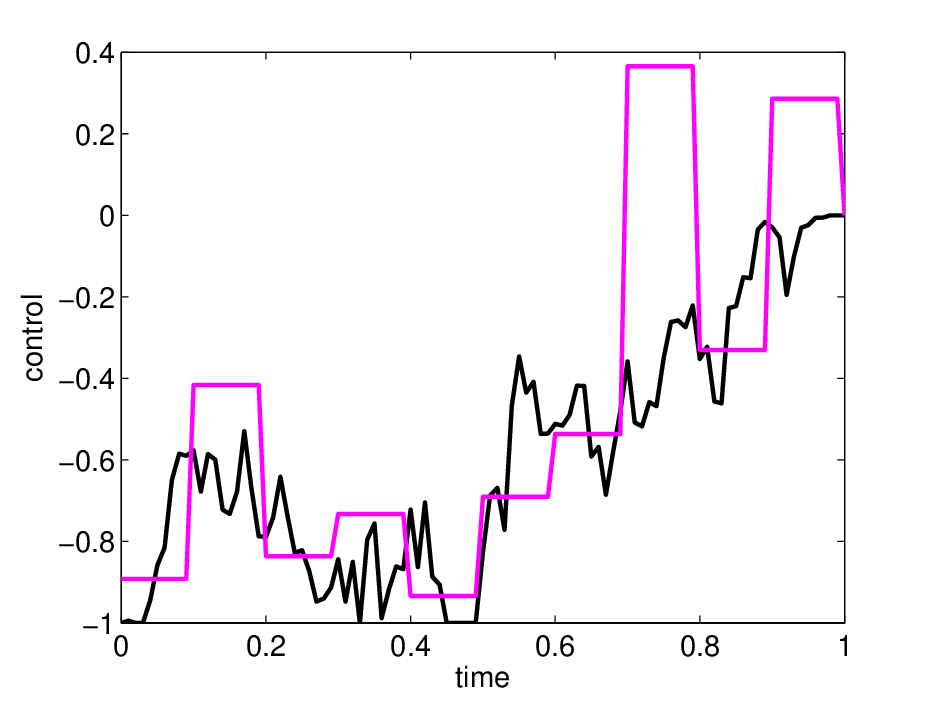}
\caption{The corresponding exact (black) optimal drift $(u(t,Y))$, given by~\eqref{eq:solution_1},~\eqref{eq:solution_2}, and approximate (magenta) optimal drift $(\widehat U_t)$, as obtained by Algorithm~\ref{alg:monte-carlo}. The approximate control is piecewise constant with interval size $\Delta = 0.1$. We see how the approximation of the optimal drift may deviate from the exact solution due to discretization and statistical variation in the Monte-Carlo sampling.}
\end{subfigure}
\caption{Minimization of the maximum of a Brownian motion with controlled drift, subject to squared control costs.}
\label{fig:minimizing_maximum}
\end{center}
\end{figure}

Algorithm~\ref{alg:monte-carlo} provides pseudocode for the approximate computation of the optimal drift $(U_t)$, employing a Monte-Carlo approximation of the expression~\eqref{eq:control_as_derivative} for $(U_t)$.\footnote{We wish to acknowledge drs. S. Thijssen (Radboud Universiteit Nijmegen) for providing us with this pseudocode and an implementation of the algorithm.} See \cite{ThijssenKappen2014} for recent results in this direction.
Parameters of the algorithm are:
\begin{itemize}
\item $\delta$: small time step for numerical simulation;
\item $\Delta$: small time step determining the interval length over which the approximate controls will be constant; $\Delta$ should be an integer multiple of $\delta$;
\item $N$: the number of samples paths generated at every iteration. 
\end{itemize}
In the algorithm, $\mathcal{N}(\mu, \sigma^2)$ denotes a (pseudo) random number generator that draws from a normal distribution with mean $\mu$ and variance $\sigma^2$.
In the experiment we set $T = 1, N = 10^2, \Delta = 0.1, \delta = 0.01$. 
An illustration of the exact and approximate computational methods of this section is provided in Figure~\ref{fig:minimizing_maximum}.

\begin{algorithm}
\caption{Monte-Carlo algorithm for relative entropy weighted minimization of the maximum of a Wiener process over a finite time interval} 
\label{alg:monte-carlo}
\begin{algorithmic}
	\State$X_0\gets 0$, $Y_0\gets 0$, $M_0\gets 0$
	\For {$t=0:\Delta:T - \Delta$}
		\For {$n = 1, 2, \ldots, N$}  \hfill Perform $N$ `rollouts' until end time $T$
			\State $X^n_t \gets Y_t$, $M^n_t \gets M_t$
			\For {$s = t, t + \delta, \ldots, T - \delta$}
				\State$X^n_{s + \delta}\gets X^n_s + \mathcal{N}(0,\delta)$
				\State$M^n_{s + \delta}\gets \max(M^n_s, X^n_{s + \delta})$ \hfill Compute the value of the cost for every rollout
			\EndFor
		\EndFor
		\State $\widehat U_t \gets \frac{\sum_{n=1}^N (X^n_{t + \Delta} - X^n_t)e^{-M^n_T}}{\sum_{n=1}^N e^{-M^n_T}}$ \hfill By law of large numbers, this is an approximation to~\eqref{eq:control_as_derivative}
		\For {$s = t, t + \delta, \ldots, t + \Delta - \delta$} \hfill Simulate a time step $\Delta$ using control $\widehat U_t$
			\State$X_{s + \delta}  \gets X_s + \mathcal{N}(0, \delta)$ \hfill $X$ is uncontrolled Wiener process
			\State$Y_{s + \delta}\gets Y_s +  \widehat U_t \delta + X_{s + \delta} - X_{s}$ \hfill $Y$ is (approximately) optimally controlled process
			\State$M_s\gets \max(M_s, Y_{s + \delta})$
		\EndFor
	\EndFor 
\end{algorithmic}
\end{algorithm}

\appendix

\section{Appendix}

\subsection{Proof of Proposition~\ref{prop:general_optimization}}

Let $(\Omega, \mathcal F, \P)$ be a probability space, and $C$ be a $\mathcal F$-measurable random variable assuming its values within the extended reals. Write $\mathcal P$ for the set of probability measures on $(\Omega, \mathcal F)$, $\mathcal P_0 := \left\{ \Q \in \mathcal P : \mbox{$\E^{\Q} [(C)^+] < \infty$ and $\mathcal H(\Q;\P) < \infty$}\right\}$.




\begin{proof}[Proof of Proposition~\ref{prop:general_optimization}]
We will distinguish two cases.
\begin{itemize}
 \item[(i)] $\E^{\P}[(C)^-] = \infty$.
In this case, essentially just using Jensen's inequality,
\begin{align*} \E^{\P}[\exp(-C)] & \geq \E^{\P} [\exp((C)^-) \1_{\{C \leq 0\}}] = \E^{\P} [\exp((C)^-) (1 -  \1_{\{C > 0\}} ) ] \\
& \geq \E^{\P} [\exp((C)^-)] - \P(C > 0) \geq \exp(\E^{\P}[(C)^-]) - \P(C > 0) = \infty.
\end{align*}
Therefore the lower bound,
\[ \inf_{\Q \in \mathcal P_0} J(\Q) \geq - \log \E^{\P} [\exp(-C)] \]
is immediate.

In the simple case $\P(C = -\infty) > 0$, we may choose $d\Q/d \P = \1_{\{C =-\infty\}}$, which establishes the required upper bound. 

Assume now $\P(C = -\infty) = 0$ and (still) $\E^{\P}[(C)^-] = \infty$.
Define measures $\Q_n$ with densities $d \Q_n/d \P = \1_{[-n,0]}(C) / \P(C \in [-n,0])$. Since $\P(C \in [-n,0]) \leq \P(C \in (-\infty, 0]))$, we have by monotone convergence
\begin{align*}
 \E^{\Q_n} \left[ C \right] & = \frac{\E \left[ C \1_{[-n,0]}(C) \right]}{\P(C \in [-n,0])} \leq \frac{\E \left[ C \1_{[-n,0]}(C) \right]}{\P(C \in (-\infty,0])} \downarrow - \frac{\E [(C)^-]}{\P(C \in (-\infty, 0])} = - \infty
\end{align*}
Furthermore $\mathcal H(\Q_n;\P) \downarrow - \log \P ( C \in [0,\infty)$ as $n \rightarrow \infty$. Therefore $J(\Q_n) \downarrow - \infty$. This minimizing sequence establishes the required upper bound.

%

\item[(ii)] $\E[(C)^-] < \infty$: Define $C_n := C \vee (-n)$. Since $C_n$ is bounded from below, $\E^{\P} \left[ \exp(-C_n) |C_n| \right] < \infty$, and $\P(C_n < + \infty) = \P(C < + \infty) > 0$, so~\eqref{hyp:finite_entropy} holds for $C_n$. Define $\mathcal P_0(C_n) = \left\{ \Q \in \mathcal P : \mbox{$\E^{\Q} (C_n)^+ < \infty$ and $\mathcal H(\Q;\P) < \infty$} \right\}$. 
Using Proposition~\ref{prop:RE_optimization},
\[ \inf_{\Q \in \mathcal P_0(C_n)} J(\Q) = - \log \E^{\P} \exp(-C_n).\]
We have trivially $\mathcal P_0 = \mathcal P_0(C_n)$. Furtermore by monotone convergence $-\log \E^{\P} \exp(-C_n) \downarrow - \log \E^{\P} \exp(-C)$ and, for $\Q \in \mathcal P_0$, $\E^{\Q} C_n \downarrow \E^{\Q} C$. Therefore, exchanging two infima,
\begin{align*} \inf_{\Q \in \mathcal P_0} J(\Q) & = \inf_{\Q \in \mathcal P_0(C_n)} \left(\inf_{n \in \N} \E^{\Q} C_n + \mathcal H(\Q;\P) \right)= \inf_{n \in \N} \left(- \log \E^{\P} \exp(-C_n) \right)\\
& = - \log \E^{\P} \exp(-C).\end{align*}
\end{itemize}

\end{proof}

\subsection*{Acknowledgements}
We are grateful to dr. ir. O. van Gaans (Mathematical Institute, Leiden University), for his comments on our work, and drs. S. Thijssen (Radboud University, Nijmegen) for providing us with his implementation of the Monte-Carlo scheme. We also wish to acknowledge dr. J. Lehec (Universit\'e Paris Dauphine) for providing us with a some details concerning his work \cite{Lehec2011}. Also we wish to acknowledge the important suggestions of the associate editor and several anonymous reviewers. These comments have had significant impact upon this paper, for which we are very much indebted.


\begin{thebibliography}{RW94b}

\bibitem[BD98]{BoueDupuis1998}
M~Bou\'{e} and P~Dupuis.
\newblock {A variational representation for certain functionals of Brownian
  motion}.
\newblock {\em The Annals of Probability}, pages 1641--1659, 1998.

\bibitem[BH13]{Banisch2013}
Ralf Banisch and Carsten Hartmann.
\newblock {Meshless discretization of LQ-type stochastic control problems}.
\newblock {\em http://arxiv.org/abs/1309.7497}, September 2013.

\bibitem[DE97]{DupuisEllis1997}
Paul Dupuis and Richard~S. Ellis.
\newblock {\em {A Weak Convergence Approach to the Theory of Large Deviations
  (Wiley Series in Probability and Statistics)}}.
\newblock Wiley-Interscience, 1997.

\bibitem[F85]{Follmer1985}
H~F\"{o}llmer.
\newblock {An entropy approach to the time reversal of diffusion processes}.
\newblock {\em Stochastic Differential Systems Filtering and Control}, 8092,
  1985.

\bibitem[F86]{Follmer1986}
H~F\"{o}llmer.
\newblock {\em {Time reversal on Wiener space}}, volume 1158 of {\em Lecture
  Notes in Math.}
\newblock Springer, Berlin, 1986.

\bibitem[Fle82]{Fleming1982}
Wendell~H Fleming.
\newblock {Logarithmic transformations and stochastic control}.
\newblock In {\em Advances in filtering and optimal stochastic control
  (Cocoyoc, 1982)}, volume~42 of {\em Lecture Notes in Control and Inform.
  Sci.}, pages 131--141. Springer, Berlin, 1982.

\bibitem[FS09]{FlemingSoner2009}
Wendell~H. Fleming and Halil~Mete Soner.
\newblock {\em {Controlled Markov Processes and Viscosity Solutions (Stochastic
  Modelling and Applied Probability)}}.
\newblock Springer, 2009.

\bibitem[HS91]{HeinricherStockbridge1991}
AC~Heinricher and RH~Stockbridge.
\newblock {Optimal control of the running max}.
\newblock {\em SIAM journal on control and optimization}, 29(4):936--953, 1991.

\bibitem[Kap05]{Kappen2005}
Hilbert~J. Kappen.
\newblock {Linear Theory for Control of Nonlinear Stochastic Systems}.
\newblock {\em Physical Review Letters}, 95(20):200201, November 2005.

\bibitem[KS91]{KaratzasShreve1991}
Ioannis Karatzas and Steven Shreve.
\newblock {\em {Brownian Motion and Stochastic Calculus (Graduate Texts in
  Mathematics)}}.
\newblock Springer, 1991.

\bibitem[Leh13]{Lehec2011}
J~Lehec.
\newblock {Representation formula for the entropy and functional inequalities}.
\newblock {\em Annales de l'Institut Henri Poincar\'{e}, Probabilit\'{e}s et
  Statistiques}, pages 1--18, 2013.

\bibitem[Mon13]{Monoyios2013}
Michael Monoyios.
\newblock {Malliavin calculus method for asymptotic expansion of dual control
  problems}.
\newblock {\em SIAM J. Financial Math.}, 4(1):884--915, 2013.

\bibitem[Nua06]{Nualart2006}
David Nualart.
\newblock {\em {The Malliavin Calculus and Related Topics}}.
\newblock 2006.

\bibitem[OK91]{OconeKaratzas1991}
Daniel~L Ocone and Ioannis Karatzas.
\newblock {A generalized clark representation formula, with application to
  optimal portfolios}.
\newblock {\em Stochastics and Stochastic Reports}, 34(3-4):187--220, 1991.

\bibitem[RW94a]{RogersWilliams1994a}
L~C~G Rogers and David Williams.
\newblock {\em {Diffusions, Markov processes, and martingales. Vol. 1}}.
\newblock Wiley Series in Probability and Mathematical Statistics: Probability
  and Mathematical Statistics. John Wiley \& Sons Ltd., Chichester, second
  edition, 1994.

\bibitem[RW94b]{RogersWilliams1994b}
L~C~G Rogers and David Williams.
\newblock {\em {Diffusions, Markov processes, and martingales. Vol. 2}}.
\newblock Wiley Series in Probability and Mathematical Statistics: Probability
  and Mathematical Statistics. John Wiley \& Sons Ltd., New York, 2nd edition,
  1994.

\bibitem[TK14]{ThijssenKappen2014}
Sep Thijssen and H.~J. Kappen.
\newblock {Path Integral Control and State Dependent Feedback}.
\newblock {\em http://arxiv.org/abs/1406.4026}, June 2014.

\bibitem[Tod06]{Todorov2006}
Emanuel Todorov.
\newblock {Linearly-solvable Markov decision problems}.
\newblock In B~Sch\"{o}lkopf, J~Platt, and T~Hoffman, editors, {\em Advances in
  Neural Information Processing Systems 19}, number~1, pages 1369--1376. MIT
  Press, Cambridge, MA, 2006.

\bibitem[Whi82]{Whittle1982}
Peter Whittle.
\newblock {\em {Optimization over time}}.
\newblock Wiley Series in Probability and Mathematical Statistics: Applied
  Probability and Statistics. John Wiley \& Sons, Ltd., Chichester, 1982.

\end{thebibliography}
\end{document}